\documentclass[11pt,bezier,epsf]{amsart}
\usepackage{amsmath,amssymb,amsfonts}
\usepackage{graphicx}
\usepackage{cite}
\usepackage{subfigure}
\newbox\subfigbox
\makeatletter
{\def\caption##1{\gdef\subcapsave{\relax##1}}%
  \let\subcapsave\@empty%
  \setbox\subfigbox\hbox%
  \bgroup}%
{\egroup%
  \subfigure[\subcapsave]{\box\subfigbox}}%
\makeatother

 \newtheorem{thm}{Theorem}[section]
 \newtheorem{cor}[thm]{Corollary}

\theoremstyle{definition}

 \newtheorem{defn}[thm]{Definition}
\theoremstyle{remark}
 
\numberwithin{equation}{section}

\begin{document}

\title[An ultraproduct method]{An ultraproduct method via left reversible semigroups to study Bruck's generalized conjecture}

\author[F. Naderi]{Fouad Naderi}

\address{Department of Mathematical and Statistical Sciences, University of
Alberta, Edmonton, Alberta, T6G 2G1, Canada.}
\email{naderi@ualberta.ca}

\subjclass[2010]{46L05; 47H10 }
\keywords{Nearly uniformly convex Banach space, Non-expansive mapping, Scattered C*-algebra, Weak fixed point property.}

\begin{abstract}
We use a method similar to ultraproducts to study the common fixed point of a
left reversible semitopological semigroup acting on a Banach space. As an application, we prove a Bruck-type theorem
for nearly uniformly convex Banach spaces to the effect that such spaces have weak fixed point property for left reversible semigroups. 
\end{abstract}

\maketitle
\section{Introduction}
Let $K$ be a subset of a Banach space $E$. A self mapping $T$ on $K$ is said to
be {\it non-expansive} if  $\|T(x)-T(y)\|\leq \|x-y\|$ for all $x,y\in K$. We
say that $E$ has the {\it weak fixed point property (weak
fpp)} if for every weakly compact convex non-empty subset $K$ of $E$, any
non-expansive self mapping on $K$ has a fixed point.

Let $S$ be a {\it semi-topological semigroup}, i.e., $S$ is a semigroup with a
Hausdorff topology such that for each $a \in S$, the mappings $s\mapsto sa$ and
$s\mapsto as$ from $S$ into $S$ are continuous. $S$ is called {\it left
reversible} if any two closed right ideals of $S$ have non-void
intersection.

An {\it action} of $S$ on a subset $K$ of a topological space $E$ is a mapping
$(s,x)\mapsto s(x)$ from $ S \times K$ into $K$ such that $(st)(x)=s(t(x))$ for
$ s,t\in S, x\in K$. The action is {\it separately continuous} if it is
continuous in each variable when the other is kept fixed. We say that
$S$ has a {\it common fixed point} in $K$ if there exists a point $x$ in $K$
such that $sx=x$ for all $s\in S$. When $E$ is a normed space, the action of $S$
on $K$ is {\it non-expansive} if $ \| s(x)-s(y)\| \leq \| x-y\|$ for all  $s \in
S$ and $x,y \in K$. There are also other types of action for a semi-topological
semigroup (see \cite{Amini} and \cite{Holmes}).

We say that a Banach space $E$ has the
{\it weak fpp for left reversible semigroups} if for every
weakly compact convex non-empty subset $K$ of $E$, any non-expansively
separately continuous action of a semi-topological semigroup $S$  on $K$ has a
fixed point.

One of the celebrated results in fixed point theory is due to Bruck \cite{Bruck}.
He has shown that if a Banach space $E$
has weak fpp, then it has weak fpp for abelian semigroups.
Now, we call the following statement {\it Bruck's Generalized Conjecture (BGC)}:

(BGC) If a Banach space $E$ has weak fpp, then it has weak fpp for any left
reversible semi-topological semigroup $S$.

In \cite{Naderi}, it has been shown that BGC is true for the preduals of von Neumann algebras.
The main results of the current work is to show BGC is also true for nearly uniformly 
convex Banach spaces.
\section{ Weak fixed point property of Bruck type}

 A Banach space $E$ is called
{\it nearly uniformly convex (NUC),} if for each $\varepsilon>0$ there exists a $\delta>0$ such that for
every sequence $(x_n)$ in the closed unit ball of $E$ with $sep(x_n):=\inf\{\rVert\ {x_n} -{ x_m}\rVert: n\neq m\}>\epsilon$ the distance $dist(0, co\{x_n\})$ is strictly less than $1-\delta$ where $co\{x_n\}$ denotes the convex hull of the sequence.

When $S$ is a left reversible semigroup, we make it to a directed set by declaring: $\alpha \geq \beta$ if and only if $\alpha S\subseteq \overline{\beta S}$. Thus, we can use $S$ as an index set for nets and speak about limit and limit-supremum with respect to this net. Also, note that when $S$ acts on a weakly compact subset $K$ of Banach space $E$, each $\alpha$ acts non-expansively on $K$ and sometimes we use the notation $T_\alpha$ instead of $\alpha$; even, we may use $\alpha.\beta$ to denote the composition $T_{\alpha}\circ T_{\beta}$.

Let $l_{\infty}(E)=\{x=(x_{\alpha}): x_{\alpha}\in E;\quad \rVert x\rVert=\sup\rVert x_\alpha \rVert <\infty \}$, and
$\mathcal{N}=\{x=(x_\alpha): x_{\alpha}\in E;\quad \lim_{\alpha} \rVert x_\alpha \rVert =0\}$. Put $\tilde{E}= l_{\infty}(E)/\mathcal{N}$ and endow it with the quotient norm $ \rVert [(x_\alpha)]\rVert=\lim \sup_{\alpha} \rVert x_\alpha \rVert$. One can embed $E$ and its subsets into $\tilde{E}$ by using constant classes. For example, for $x\in K$ let $\dot{x}=[(x)]$ denotes the equivalence class containing the constant net $(...,x,x,x,...)$. So, $\dot{K}=\{ \dot{x}: x\in K \} $ is a subset of $\tilde{E}$. Also, we define $\tilde{K}=\{[(k_\alpha)]\in \tilde{E}: k_{\alpha}\in K \quad \text{for each } \alpha \}$; and note that, $\dot{K}\subseteq \tilde{K}$. The process of embedding preserves the properties of being closed and convex for subsets of $E$. If each $T_\alpha$ is non-expansive on $K$, then the mega mapping
$\tilde{T}:\tilde{K}\longrightarrow \tilde{K}$ defined by $\tilde{T}[(x_\alpha)]=[(T_{\alpha}x_\alpha)]$ is also non-expansive.
 Another piece of notation, when $(x_{k})$ is a sequence in $E$, then $f_{k}=[(\alpha.x_{k})_{\alpha}]$ is a sequence in $\tilde{E}$.

\begin{defn} \label{asym}
(a) Let $K$ be a non-empty subset of a Banach space $E$ and $(x_\alpha)_{\alpha\in S}$ be a bounded net in $E$. For
each $k\in K, \alpha \in S$ define
$$r(K, (x_\alpha))=\inf\{\lim\sup_{\alpha}\parallel x_{\alpha}-k\parallel:  k\in K \}$$
$$ AC(K,(x_\alpha))=\{ k\in K : \lim\sup_{\alpha}\parallel x_{\alpha}-k\parallel= r(K, (x_\alpha)) \}$$

The set $AC(K,(x_\alpha))$ (the number
$r(K, (x_\alpha))$ ) will be called the {\it asymptotic center (asymptotic radius)} of $(x_\alpha)_{\alpha \in S}$ in $K$. 
These are the generalizations of Chebyshev center and radius and are due to Edelstein\cite{Edelstein}.
The asymptotic center is always non-empty for weakly compact set $K$.

(b) When viewed in $\tilde{E}$, these notions are seen as:

$$r(\dot{K}, [(x_\alpha)])=\inf\{\parallel [(x_{\alpha})]-\dot{k}\parallel  k\in K \},$$
$$ AC(\dot{K}, [(x_\alpha)])=\{ \dot{k}\in \dot{K} : \parallel [(x_{\alpha})]-\dot{k}\parallel=r(\dot{K}, [(x_\alpha)]) \}$$.

\end{defn}

The proof of the following theorem uses some ideas from \cite{Butsan} and \cite{Wisnicki}.
\begin{thm}\label{Bruck}
Let $K$ be a bounded closed convex non-empty subset of a Banach space $E$, and let $S$ be a left reversible semigroup acting non-expansively and separately continuous on $K$. Suppose there exists a $\lambda\in [0,1)$ such that for each net
$(x_\alpha)$ in $K$ and each net $(y_\alpha)$ in $AC(K, (x_\alpha))$ we have $r(K, (y_\alpha))\leq \lambda r(K, (x_\alpha))$. Then, $S$ has a common fixed point in $K$.
\end{thm}
{\bf Proof.} Let $x_{0} \in K$. Put $r_{0}=r(K,(\alpha.x_{0}))$ and $A_{1}=AC(K,(\alpha.x_{0}))$. Let $x_{1}\in A_{1}$. Since the action is non-expansive we see that $\alpha.x_{1}\in A_{1}$ and $[(\alpha.x_{1})]\in \tilde{A_{1}}$. Put  $r_{1}=r(K,(\alpha.x_{1}))$ and $A_{2}=AC(K,(\alpha.x_{1}))$. By induction,
we get a sequence $x_{k}$ such that $x_{k}\in A_{k}$, $[(\alpha.x_{k})]\in \tilde{A_{k}}$, $r_{k}=r(K,(\alpha.x_{k}))$, $A_{k+1}=AC(K,(\alpha.x_{k}))$ and $r_{k}\leq {\lambda}^{k}r_{0}$. As in the proof of 
\cite[Theorem 3.17]{Butsan} and \cite[Theorem 5.3]{Wisnicki}, there exists a $u\in K$ such that the sequence  $f_{k}=[(\alpha.x_{k})]$ converges to $\dot{u}$ in $\tilde{K}$. Note that the use of $\alpha\in S$ instead of $n\in \mathbb{N}$ in those theorems for the part we need, cause no problem since the convergence occurs in $k$. So, the limit of
$[(\alpha.x_{k})]$ is $[(u)]$ regardless of the indexing $\alpha$s.  Let $\beta\in S$ be arbitrary and consider the elements $[(\beta.\alpha.x_{k})]$ and the constant element $[(\beta.u)]$ which are indexed by $\alpha$s. Since the action is non-expansive
we get the $\lVert [(\beta.\alpha.x_{k})]-[(\beta.u)] \rVert \leq \lVert [(\alpha.x_{k})]-[(u)] \rVert$. Hence the sequence $[(\beta.\alpha.x_{k})]$ converges to $[(\beta.u)]$ in $\tilde{K}$.  But, $\beta.\alpha \in S$, so by the first part of the current proof, $[(\beta.\alpha.x_{k})]$ is a another sequence like $[(\alpha.x_{k})]$ with the same limit. That is, by the same discussion, $[(\beta.\alpha.x_{k})]$ must have the limit $[(\beta.u)]=[(u)]$. Hence, $(\beta.u-u)_{\alpha}\in \mathcal{N}$. Therefore, $\beta.u=u$, which is the desired result.$\blacksquare$
\begin{cor}
Let $K$ be a weakly compact convex non-empty subset of a  nearly uniformly convex Banach space $E$, and let $S$ be a left reversible semigroup acting non-expansively and separately continuous on $K$. Then, $S$ has a common fixed point in $K$.
\end{cor}
{\bf Proof.} The set $K$ is closed and bounded. Also, nearly uniformly convex Banach spaces satisfy the inequality in Theorem \ref{Bruck}. So, the result follows.$\blacksquare$


\end{document}